\RequirePackage{fix-cm}
\documentclass[smallextended]{svjour3}       
\smartqed  
\usepackage{graphicx}
\usepackage{amsmath,amssymb,amsfonts}

\newcommand{\one}[1]{\mathbf{1}_{#1}}

\newcommand{\VecSet}[2]{\mathbb{#1}^{#2}}

\newcommand{\Encode}[2]{\mathcal{C}_{#1,#2}}
\newcommand{\Simplex}[1]{\mathbb{S}^{#1}}
\newcommand{\SimMat}[1]{\mathbb{S}^{#1 \times #1}}
\newcommand{\SimMats}[2]{\mathbb{S}^{#1 \times #2}}

\newcommand{\Trans}[2]{#1^{(#2)}}
\newcommand{\TransBlock}[3]{#1_{#3}^{(#2)}}
\newcommand{\TransBlockVec}[4]{#1_{N(#3-1) + #4}^{(#2)}}

\newcommand{\zero}[1]{\mathbf{0}_{#1}}

\begin{document}

\title{Recursive Markov Process for Iterated Games with Markov Strategies
}

\titlerunning{Recursive Markov Process}        

\author{Shohei Hidaka  
}


\institute{S. Hidaka \at
              1-1 Asahidai, Nomi, Ishikawa, Japan \\
              Tel.: +81-761-51-1717\\
              Fax: +81-761-51-1775\\
              \email{shhidaka@jaist.ac.jp}           
}

\date{Received: date / Accepted: date}

\maketitle

\begin{abstract}
The dynamics in games involving multiple players, who adaptively learn from their past experience, is not yet well understood.
We analyzed a class of stochastic games with Markov strategies in which players choose their actions probabilistically. 
This class is formulated as a $k^{\text{th}}$ order Markov process, in which the probability of choice is a function of $k$ past states.
With a reasonably large $k$ or with the limit $k \to \infty$, numerical analysis of this random process is unfeasible.
This study developed a technique which gives the marginal probability of the stationary distribution of the infinite-order Markov process, which can be constructed recursively. 
We applied this technique to analyze an iterated prisoner's dilemma game with two players who learn using infinite memory.
\end{abstract}

\section{Games with adaptive learners \label{sec-IPD}}

Conventional game theory treats each player as a rational decision maker who is making decisions based on sufficient information regarding the game being played. The dynamics in such a game are generally characterized by its Nash equilibria --- states in which none of the agents may profit by changing their actions \cite{Nash1950}.

Real social problems, however, are often vastly more complex than such formulations \cite{Camerer1999,Camerer2003,Camerer2010,Roth1995}.
In reality, each agent has limited computational resources and limited information about the game it is playing. Under such constraints, learning -- making inductive inferences of future behaviour from limited experience of past behaviour -- plays a crucial role in finding locally optimal actions. One key question regarding games in which agents have limited information but employ learning techniques concerns their long-term equilibria \cite{Camerer2003}.

A class of iterative games with {\em reinforcement learning} \cite{Sutton1998} has been investigated in both theoretical  \cite{Borgers1997,Sato2002} and empirical studies \cite{Camerer1999,Camerer2003,Camerer2010,Roth1995}.
In this class of games, the only information available to each agent is a number of its own actions and the rewards of these actions. The probabilities for the players' next actions are computed according to the weighted averages of rewards for the possible actions.

As each player's decision is probabilistic, a game of this class in which players have access to their history going back $k$ turns can be formulated as a $k^{\text{th}}$-order finite Markov process. The equilibria of such a game are then characterized by its stationary distributions.
Although this formulation is mathematically simple, not many past studies have taken this approach to the analysis of iterative games with reinforcement learning. The reason for this is that the number of states in an iterated game between learners with $k$-step memory grows exponentially with $k$. To our knowledge, only a few special cases, such as those in which $k=1$ \cite{Nowak1990} and some specific games in which $k$ is large \cite{Hidaka2015}, have been analyzed.

The present study shows a new technique which approximately computes a marginal stationary distribution of a $k^{\textrm{th}}$-order Markov process with sufficiently large $k$.
This technique is applied to an iterated game of prisoners' dilemma with learning players.
Our analysis demonstrates the utility of the proposed computational technique as a numerical tool for the analysis of dynamic games with learners.
The main result of this paper is that, for a certain class of iterated game involving players having Markovian strategies with the recursive property,
its marginal stationary distribution $\omega \in \mathbb{R}^{N}$ of $k^{\text{th}}$-order Markov process approaches to a form as $k \to \infty$:
\begin{equation}
\label{eq-MainResults}
 \omega = Q( \omega ) \omega, 
\end{equation}
where the transition matrix $Q( \omega ) \in \mathbb{R}^{N \times N}$ is a function of $\omega$.
This form gives a finite form of equality for the desired marginal stationary distribution, and analysis of it would reveal the long-run nature of the given stochastic game.

In Section \ref{sec-Game}, we describe the iterated version of prisoners' dilemma and reinforcement learning as a motivating problem for development of the mathematical technique.
In Section \ref{sec-Markov}, we introduce a Markov process formulation for the iterated prisoners' dilemma.
In Section \ref{sec-Sketch}, we present the main result of the paper with regard to the marginal stationary distribution of an infinite order Markov process.
In Section \ref{sec-Recursive} we define a new class of Markov process motivated by the main result.
In Section \ref{sec-Application}, we demonstrate the new computational technique by applying it to the iterated prisoners' dilemma.

\section{Iterated prisoners' dilemma \label{sec-Game}}

The prisoners' dilemma is a classical game which has long been used as a minimal model demonstrating the difficulty of mutual cooperation. In the prisoners' dilemma, each of two players chooses an action from either Cooperation (C) or Defection (D), and each player is given a certain payoff depending on the actions of both players. Each player benefits fairly when both choose Cooperation (CC). However, one player can gain even more by choosing Defection on the condition that the other chooses Cooperation. This incentivizes Defection for both players, and the game results in mutual Defection (DD) with individual payoffs lower than those in the mutual Cooperation case. This is the only Nash equilibrium in the prisoners' dilemma. The "dilemma'' is that two rational players cannot escape from mutual Defection with unhappy payoffs, even though mutual Cooperation is the more beneficial option.

The basic game has been extended to games with multiple agents, iterated steps, stochastic strategies, situations affected by noise, and certain topologies for agent interactions \cite{Axelrod1984,Galla2009,Sandholm1996}. In the iterative variant of these models, each agent can adaptively choose its action on the basis of a series of past actions and payoffs. One of the simplest cases is completely analyzed based on the finite Markov formalism \cite{Nowak1990}, but more general cases remain open for further research.

In a recent study \cite{Torii2014}, we analyzed the iterated prisoners' dilemma (IPD) with two probabilistic learners. This analysis suggested that IPD with two learners can result in mutual cooperation, but this is limited only to learners with short memories due to the previously mentioned computational difficulty involved in analyzing iterated games involving players with long memories. The present study develops a new computational tool which extends this analysis to learners with long memories.

\subsection{IPD with reinforcement learners}
We briefly outline the IPD with learners who learn using reinforcement learning. 

\begin{definition}[Iterated prisoners' dilemma]

\label{def-LearningGames}

We label the players $1$ and $2$, and label the two moves available to each player at each step as $0$ (for Cooperation) and $1$ (for Defection). For each integer $t$ and for each $i = 1,2$, we denote by $x_{t,i}$ the choice of action that player $i$ made on turn $t$.
Define a function $f : \{0,1\}^2 \rightarrow \{1, 2, 3, 4\}$ by
\begin{equation}
\label{eq-Indexingmap}
 f(a_1,a_2) := 2a_1 + a_2 + 1.
\end{equation}
The function $f$ encodes each of the four possible outcomes on a single turn of the game as the integers $1, 2, 3,$ and $4$. We define $X_t := f(x_{t,1}, x_{t,2})$.
We write $ X_{t}^{t+k-1} := \left( X_{t}, X_{t+1}, \ldots, X_{t+k-1} \right)$
for a sequence of $k$ states from step $t$.

We assume the existence of a payoff map $r = (r_1, r_2) : \{1, 2, 3, 4\} \rightarrow \mathbb{R}^2$ such that $r_i(X_t)$ is the payoff to the $i^{\textrm{th}}$ agent resulting from the actions taken by both agents at turn $t$. Note that, under our assumption, the payoff scheme remains constant over time.
We further assume that the payoffs are symmetric so that, for any $a, b \in \{0, 1\}$,
\[r_1(f(a,b)) = r_2(f(b,a)) =: R_{a,b}.\]

An iterated prisoners' dilemma is a game which satisfies these assumptions as well as the inequalities
\[R_{01} < R_{11} < R_{00} < R_{10},\ \textrm{and}\]
\[ R_{01} + R_{10}  < 2R_{00}.\]

\end{definition}


\begin{definition}[IPD with reinforcement learners]
In the reinforcement learning model, we begin with an IPD as defined above. However, each agent chooses an action based on a function of the rewards it received for its past $k$ actions, with $k$ specified at the outset. Specifically, for agent $i$, this function is
\[
 \phi_{i,x}^{\alpha_{i}}\left( X_{t-k}^{t-1} \right)
 = \sum_{s=1}^{k}\alpha_{i}^{s} \delta_{x,x_{t-s, i}} r_{i}( X_{t-s} ),
\]
where
$\alpha_{i} \in [0, 1]$ is a memory-retention parameter, and 
$\delta_{x,y}$ is the Kronecker delta, which takes the value 1 when $x$ and $y$ agree and is 0 otherwise.

Using this weighted rewards function along with a sensitivity parameter $\beta_{i} \ge 0$,
the probability that the $i^{\text{th}}$ agent chooses action $x$ at step $t$ is
\begin{equation}
\label{eq-ChoiceIndividual}
 P\left( x \mid X_{t-k}^{t-1} \right) 
=
 \frac{ 
   \exp\left( \beta_{i} \phi_{i,x}^{\alpha_{i}}\left( X_{t-k}^{t-1} \right) \right) 
 }{
   \sum_{x = 0}^{1}\exp\left( \beta_{i} \phi_{i,x}^{\alpha_{i}}\left( X_{t-k}^{t-1} \right) \right) 
 }
\end{equation}

We assume that the agents choose their actions independently at each turn so that, for any $(a_1, a_2) \in \{0, 1\}^2$,
\begin{equation}
\label{eq-Choice}
 P\left( X_{t} = f(a_1, a_2) \mid X_{t-k}^{t-1} \right) 
=
\prod_{i = 1}^{2} P\left( a_i \mid X_{t-k}^{t-1} \right).
\end{equation}
\end{definition}


\section{Markov process \label{sec-Markov}}
Equation (\ref{eq-Choice}) allows us to calculate the conditional probabilities $P\left( X_{t-k+1}^{t} \mid X_{t-k}^{t-1} \right)$. In this manner, we construct a Markov chain with states consisting of all possible length $k$ move sequences of the players in an IPD with reinforcement learning. This is the $k^{\text{th}}$-order Markov process corresponding to the variables $X_{t}$. We encode the states in this $k^{\text{th}}$-order Markov process with the integers $1, 2, \ldots, 4^k$ and describe its transition matrix with respect to this encoding as follows.
In a previous work under the same formulation of stochastic game, Hidaka and colleagues analyzed the rock-scissor-paper game by approximating it using a finite Markov process with a small $k$ \cite{Hidaka2015}.

\subsection{Transition matrix}
Here let us consider a class of iterated games with Markov strategies, which includes the IPD introduced in the previous section as a special case. 
Let us write the set of $N$ states at each step of a Markov process by $\mathcal{N} = \{ 1, 2, \ldots, N \}$.
For each $i$, let us write the random variable indicating the state at $i^{\text{th}}$ step by $X_{i} \in \mathcal{N}$ and the series of $k$ states by 
$X = ( X_{1}, X_{2}, \ldots, X_{k} )~\in ~\mathcal{N}^{k}$. 
For the $k^{\text{th}}$ order Markov process,
we assign each series of states in $\mathcal{N}^{k}$ an integer in the set $\Encode{N}{k} := \{ 1, 2, \ldots, N^k \}$
by the indexing map $h_{N,k}: \mathcal{N}^{k} \mapsto \Encode{N}{k}$ defined by
\begin{equation}
\label{eq-IndexingMap} \nonumber
 h_{N,k}\left( X \right) := 
 1 + \sum_{j = 1}^{k}( X_{j}-1 )N^{j}.
\end{equation}
For a $k^{\text{th}}$ order Markov process, we call each integer $i \in \Encode{N}{k}$ {\em code} and the map $h_{N,k}$ {\em encoder}.
Without loss of generality, we keep this encoding system for the set of series of states throughout this paper.

Let us write a series of states from $t$ to $s$ by $X_{t}^{s} := \left( X_{t}, X_{t+1}, \ldots, X_{s} \right)$.
For $i \in \Encode{N}{k}$, denote the set of integers, which are codes for series of states $X_{1}^{k}$ transitable from $X_{0}^{k-1} = h_{N,k}^{-1}\left( i \right)$, by
\[
 \mathcal{H}_{N,k}(i) := \left\{ h_{N,k}\left( X_{1}^{k} \right) : i = h_{N,k}\left( X_{0}^{k-1} \right)\right\}.
\]
For each $i$, the set $\mathcal{H}_{N,k}(i)$ 
consists of the $N$ codes for the $k^{\textrm{th}}$-order Markov process. 
First we define 
the transition matrix $\Trans{Q}{k} \in \SimMat{N^k}$ for the $k^{\textrm{th}}$-order Markov process with respect to this encoding as follows.
\begin{definition}[Transition matrix]
\label{def-transition} 
The transition matrix $\Trans{Q}{k} \in \SimMat{N^k}$ for the $k^{\text{th}}$-order Markov process is defined by
\begin{equation}
\label{eq-TransitionMatrix} \nonumber
 \left( \Trans{Q}{k} \right)_{i,j} := P\left( h_{N,k}^{-1}(i) \mid h_{N,k}^{-1}(j) \right),
\end{equation}
where $\left( Q \right)_{i,j}$ is the $(i, j)$ element of the matrix $Q$.
Observe that, unless $i \in \mathcal{H}_{N,k}(j)$, $\left( \Trans{Q}{k} \right)_{i,j} = 0$.
\end{definition}

In order to analyze the properties of the transition matrix $\Trans{Q}{k}$, 
let us introduce notations for the vectors and matrices as follows.
Let $\Simplex{N} := \{ (x_{1}, x_{2}, \ldots, x_{N} )^{T}: x_{i} \ge 0 \text{ for every } i, \text{ let } \sum_{i}^{N}x_{i} = 1\}$ be $N-1$ dimensional simplex,
and let $\SimMat{N} := \left( y_{1}, y_{2}, \ldots, y_{N} \right)$ denote a simplex matrix with simplex vectors $y_{i} \in \Simplex{N}$. 
Let us denote the zero vector by $\zero{N} = (0, 0, \hdots, 0)^{T} \in \VecSet{R}{N}$, 
the identity matrix by $E_{N} \in \SimMat{N}$ and the unit vector by
\[
 e_{N,i} := \left( 0,\hdots,0,\stackrel{i}{\stackrel{\vee}{1}},0,\hdots,0 \right)^{T} \in \Simplex{N},
\]
and let
$E_{N,i} := e_{N,i} e_{N,i}^{T}$.
We define a special permutation matrix called the {\em commutation matrix} \cite{Magnus1988b} by:
\begin{equation}
\label{eq-Commutation} \nonumber
 C_{n, m} := \sum_{i = 1}^{m} e_{m, i}^{T} \otimes E_{n} \otimes e_{m, i}.
\end{equation}
where $\otimes$ denotes the Kronecker product.

For the state $i \in \Encode{N}{k}$ of a $k^{\text{th}}$ order Markov process with $\Trans{Q}{k} \in \SimMat{N^k}$ 
and its corresponding indices $m_{1} < \ldots < m_{N} \in \mathcal{H}_{N,k}(i)$,
we define 
\[
 \TransBlock{q}{k}{i} := \left( P( m_1 \mid i ), \ldots, P( m_N \mid i ) \right)^{T}\in \Simplex{N},
\]
and for $i \in \Encode{N}{k-1}$, we define the simplex matrix
\[
 \TransBlock{Q}{k}{i} := \left( \TransBlock{q}{k}{N(i-1)+1}, \ldots, \TransBlock{q}{k}{N(i-1)+N} \right) \in \SimMat{N}.
\]
Using the above notation, Hidaka and colleagues \cite{Hidaka2015Simple} showed that
an arbitrary $k^{\text{th}}$ order transition matrix $\Trans{Q}{k}$ can be decomposed as:
\begin{equation}
 \label{eq-TransitionMat}
  \Trans{Q}{k} = C_{N,N^{k-1}} \sum_{i \in \Encode{N}{k-1}} E_{N^{k-1},i} \otimes \TransBlock{Q}{k}{i}.
\end{equation}
To illustrate this specifically, 
consider a small example with $N=2$,  $k=2$, 

\[
Q_{1}^{(2)} :=  \begin{pmatrix}
q_{0, 0 \mid 0, 0} & q_{1, 0 \mid 0, 1} \\
q_{0, 1 \mid 0, 0} & q_{1, 1 \mid 0, 1}  
\end{pmatrix} 
\text{ and } 
Q_{2}^{(2)} := 
\begin{pmatrix}
q_{0, 0 \mid 1, 0} & q_{1, 0 \mid 1, 1} \\
q_{0, 1 \mid 1, 0} & q_{1, 1 \mid 1, 1}  
\end{pmatrix},
\]
 where each non-zero element 
$q_{t+1, t+2 \mid t, t+1}$ is the transition probability from $(X_{t}, X_{t+1}) \in \{0, 1\}^{2}$ to the state $(X_{t+1}, X_{t+2}) \in \{0, 1\}^{2}$. Then it has the decomposition of transition matrix 
\[
Q^{(2)} = 
\overbrace{ \begin{pmatrix}
1 & 0 & 1 & 0 \\
1 & 0 & 1 & 0 \\
0 & 1 & 0 & 1 \\
0 & 1 & 0 & 1
\end{pmatrix} 
}^{ C_{2,2} }
\overbrace{
\left(
\begin{array}{cc}
Q_{1}^{(2)}
& \text{\Large{\:0}} \\
\text{\Large{\:0}} 
& 
Q_{2}^{(2)}
\end{array}
\right)
}^{ \sum_{i=1,2}E_{2,i} \otimes Q_{i}^{(2)} }
= 
\left(
\begin{array}{cc}
 \begin{matrix}
q_{0, 0 \mid 0, 0} & 0 & q_{0, 0 \mid 1, 0} & 0 \\
q_{0, 1 \mid 0, 0} & 0 & q_{0, 1 \mid 1, 0} & 0 \\
0 & q_{1, 0 \mid 0, 1} &  & q_{1, 0 \mid 1, 1} \\
0 & q_{1, 1 \mid 0, 1}   & 0 & q_{1, 1 \mid 1, 1} 
\end{matrix} 
\end{array}
\right).
\]

\section{Analysis of a higher order Markov process \label{sec-Sketch}}
\subsection{Sketch}
Our goal is to develop a technique with which we can analyze the marginalized stationary probability of the $K^{\text{th}}$ Markov process in the limit $K \to \infty$ or a sufficiently large but finite $K$.
That is, more formally, $M \theta_{K}$ ($M \in \mathbb{R}^{N \times N^{K}}$ is an appropriate marginalization matrix, defined shortly), where $\theta_{K} = Q^{(K)} \theta_{K}$ is the stationary vector for a given transition matrix $Q^{(K)}$.
As it has no expression in a closed form in general, we will explore a condition under which such marginal stationary distribution $M \theta_{K}$ can be expressed in a closed form. 
As multiple argument steps are required to understand this characteristic of the Markov process, it would be useful to have a schematic sketch here.

 The key observation comes from the decomposition of transition matrix (\ref{eq-TransitionMat}). 
In the following section, we extend this idea by expressing three types of linear operators, {\em branching, marginalization,} and {\em cycling}, each represented in a matrix form. The decomposition of transition matrix (\ref{eq-TransitionMat}) is identical to a multiplication of these three types of matrices. 

We will extend the transition matrix by defining a {\em $k$-shift} matrix. A transition matrix (\ref{eq-TransitionMat}) is identified as a $1$-shift matrix, and the $k$-shift matrix is expressed by multiplication of a cycling matrix and a $k$ series of the branching and marginalization matrices. In general, a $k$-shift matrix is a linear operator which shifts a marginal probability 
$P(X_{t})$ to $P(X_{t+k})$ under a given $K^{\text{th}}$ order Markov process. 
Then, we consider ``stationary vector of the shift matrix'' $\omega_{k} \in \Simplex{N}$ for $k$-shift matrix $S_{k} \in \SimMat{N}$, which satisfies $\omega_{k} = S_{k} \omega_{k}$. 

The reader should be aware that $\omega_{k}$ is a ``stationary vector (the primary eigenvector) of $k$-shift matrix'', which may not be identical to the ``marginalized stationary probability of $K^{\text{th}}$ Markov process'' $M \theta_{K}$ above.
In a later section, we will show that these two types of vectors, however, essentially have a similar recursive property in common. 
Thus, our main question is what condition makes the stationary of marginal probability $\omega_{k}$ equivalent to the marginalized stationary probability $M \theta_{K}$, if possible.

In summary, our argument takes the steps as follows:
\begin{enumerate}
 \item Define branching, marginalization, and cycling matrices (Section \ref{sec-operators}).
 \item Define the $k$-shift matrix and show the recursive property of the stationary vector for it (Section \ref{sec-shift} and Lemma \ref{prop-RecursiveForm}).
 \item Show the recursive property of the stationary vector of the $K^{\text{th}}$ order Markov process (Lemma \ref{prop-marginal}).
 \item Show a condition under which the stationary vector of the shift matrix is identical to the marginalized stationary vector (Theorem \ref{thm-RecursiveMarginal}).
\end{enumerate}
The result obtained in this section motivates us to define the class of recursive Markov process. This process can be constructed by a recursive procedure and will be defined in Section \ref{sec-Recursive}. 

\subsection{Shift matrix and linear operators \label{sec-operators}}
In this section, we introduce three linear operators in their matrix forms to show the basic properties of the $k^{\text{th}}$ order transition matrix
serving as the foundation of our main result.

\begin{definition}
For $0 \le m \le k$, we define the $k^{\text{th}}$ order {\em marginalization matrix} 
\[
  \TransBlock{M}{k}{m} :=  E_{N^{m}} \otimes \one{N}^{T} \otimes E_{N^{k-m}} 
\in \SimMats{N^{k}}{N^{k+1}}.
\]
For $0 \le m \le k$ and 
the tuple of vectors $\TransBlock{\mathcal{Q}}{k}{N} := \left( \TransBlock{q}{k}{1}, \ldots, \TransBlock{q}{k}{N^k} \right)$,
we define the $k^{\text{th}}$ order {\em branching matrix} 
\[
 \TransBlock{B}{k}{m}\left( \TransBlock{\mathcal{Q}}{k}{N} \right) 
  := \sum_{i \in \Encode{N}{m}, j \in \Encode{N}{k-m}} E_{N^{m}, i } 
          \otimes \TransBlock{q}{k}{ N^{k-m} ( i - 1 ) + j } \otimes E_{N^{k-m}, j }
 \in \SimMats{N^{k+1}}{N^{k}}.
\]
For $0 \le m \le k$, we define the {\em cycling matrix}
\[
\TransBlock{C}{k}{m} := C_{N^{m}, N^{k-m}} = \sum_{i \in \Encode{N}{k-m}} e_{N^{k-m}, i}^{T} \otimes E_{N^{m}} \otimes e_{N^{k-m}, i}. 
\] 
\end{definition}

For example with $N = 2$, we have a marginalization matrix and a cycling matrix
\begin{equation}
\TransBlock{M}{2}{2}
= 
\begin{pmatrix}
\one{2}^{T} & & &
\\
& \one{2}^{T} & &
\\
& & \one{2}^{T} & &
\\
& & & \one{2}^{T}
\end{pmatrix} 
\text{ and }
\TransBlock{C}{2}{1}
= 
\begin{pmatrix}
1 & 0 & 0 & 0
\\
0 & 0 & 1 & 0
\\
0 & 1 & 0 & 0
\\
0 & 0 & 0 & 1
\end{pmatrix} 
\end{equation}
and a branching matrix for the tuple $\TransBlock{\mathcal{Q}}{2}{2} = \left( \TransBlock{q}{2}{1}, \TransBlock{q}{2}{2}, \TransBlock{q}{2}{3}, \TransBlock{q}{2}{4} \right)$
\begin{equation}
\TransBlock{B}{2}{2}\left( \TransBlock{\mathcal{Q}}{2}{2} \right)
= 
\left(
\begin{array}{cc}
\begin{matrix}
\TransBlock{q}{2}{1} &  \\
 & \TransBlock{q}{2}{2}
\end{matrix}
& \text{\Large{\:0}} \\
\text{\Large{\:0}} 
& 
\begin{matrix}
\TransBlock{q}{2}{3} &  \\
 & \TransBlock{q}{2}{4}
\end{matrix}
\end{array}
\right).
\end{equation}

Let each element of the vector $\Trans{\theta}{k} \in \Simplex{N^k}$ be
an arbitrary stochastic vector consisting of the probability $P\left( X_{1}, X_{2}, \ldots, X_{k} \right) $ for
$h_{N,k}\left( X_{1}, X_{2}, \ldots, X_{k} \right) \in \Encode{N}{k}$.
Let $\TransBlock{\mathcal{Q}}{k}{N} := \left( \TransBlock{q}{k}{1}, \ldots, \TransBlock{q}{k}{N^k} \right)$
be the tuple of simplex vectors, such that the vector
$\left( \left( \TransBlock{q}{k}{1} \right)^{T}, \ldots, \left( \TransBlock{q}{k}{N^k} \right)^{T} \right)^{T}$ consists of 
the conditional probability $P\left( Y \mid X_{1}, X_{2}, \ldots, X_{k} \right)$.
Then, the three types of matrices introduced above correspond to the operators on the stochastic vector $\Trans{\theta}{k}$ as follows.
\begin{enumerate}
 \item{
  Marginalization $\TransBlock{M}{k-1}{m}\Trans{\theta}{k} \in \Simplex{N^{k-1}}$: 
  $P\left( X_{1}, \ldots, X_{m}, X_{m+2}, \ldots, X_{k} \right)$.
  }
 \item {
  Branching $\TransBlock{B}{k}{m}\left( \TransBlock{\mathcal{Q}}{k}{N} \right) \Trans{\theta}{k} \in \Simplex{N^{k+1}}$: $P\left( X_{1}, \ldots, X_{m}, Y, X_{m+1}, \ldots, X_{k} \right) $
 }
 \item{
  Cycling $\TransBlock{C}{k}{m} \Trans{\theta}{k} \in \Simplex{N^k}$: 
$P\left( X_{m+1}, X_{m+2}, \ldots, X_{k}, X_{1}, \ldots, X_{m} \right) $
 }
\end{enumerate}

\begin{figure}[h]
\begin{center}
\includegraphics[width=10cm]{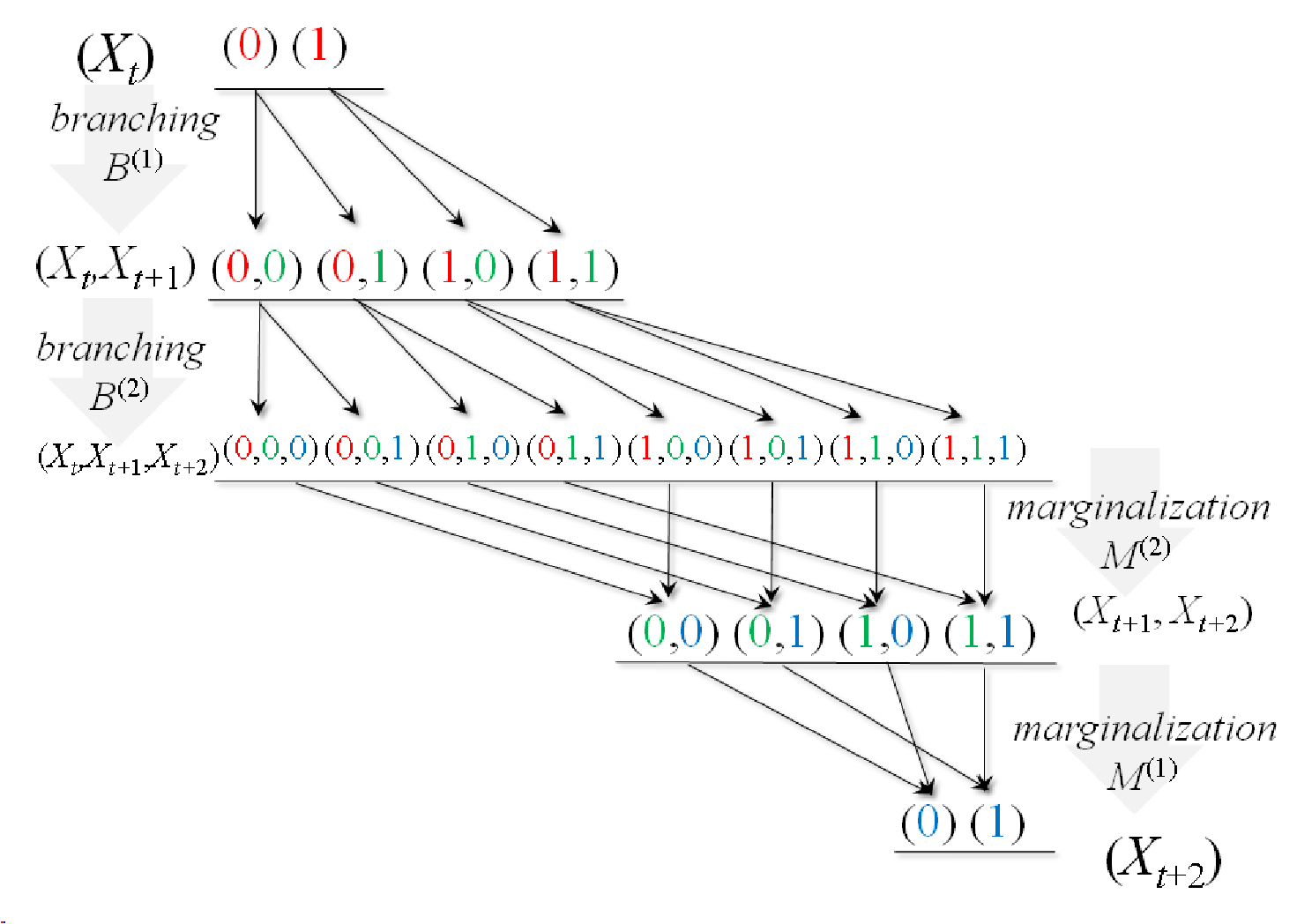}
\caption{Branching and marginalization matrix for $N=2$ and $k=1,2,3$.}
\label{fig-BandM}
\end{center}
\end{figure}
Figure \ref{fig-BandM} illustrates the branching and marginalization matrices
for $N=2$ and $k=1, 2$.
Consider, for example, the tuples for $k=1, 2$:
\[
\TransBlock{\mathcal{Q}}{1}{2} := 
\left(
\overbrace{ 
 \begin{pmatrix}
q_{0, 0 \mid 0} \\
q_{0, 1 \mid 0}
\end{pmatrix} 
}^{ \TransBlock{q}{1}{1} }
,
\overbrace{
 \begin{pmatrix}
q_{1, 0 \mid 1} \\
q_{1, 1 \mid 1}  
\end{pmatrix} 
}^{ \TransBlock{q}{1}{2} }
\right)
\],
and 
\[
\TransBlock{\mathcal{Q}}{2}{2} := 
\left(
\overbrace{
 \begin{pmatrix}
q_{0, 0, 0 \mid 0, 0} \\
q_{0, 0, 1 \mid 0, 0}
\end{pmatrix}
}^{ \TransBlock{q}{2}{1} }
,
\overbrace{
 \begin{pmatrix}
q_{0, 1, 0 \mid 0, 1} \\
q_{0, 1, 1 \mid 0, 1}  
\end{pmatrix} 
}^{ \TransBlock{q}{2}{2} }
,
\overbrace{
 \begin{pmatrix}
q_{1, 0, 0 \mid 1, 0} \\
q_{1, 0, 1 \mid 1, 0}
\end{pmatrix}
}^{ \TransBlock{q}{2}{3} }
,
\overbrace{
 \begin{pmatrix}
q_{1, 1, 0 \mid 1, 1} \\
q_{1, 1, 1 \mid 1, 1}  
\end{pmatrix} 
}^{ \TransBlock{q}{2}{4} }
\right)
\]
 where each non-zero element 
$q_{ X_t, X_{ t+1 } \mid X_{ t } }$ (and $q_{ X_{ t }, X_{ t+1 }, X_{ t+2 } \mid X_{ t }, X_{ t+1 } }$) is the transition probability from $(X_{t}) \in \{0, 1\}^{2}$ to the state $(X_{t}, X_{t+1}) \in \{0, 1\}^{2}$ (and from $(X_{t}, X_{t+1}) \in \{0, 1\}^{2}$ to the state $(X_{t}, X_{t+1}, X_{t+2}) \in \{0, 1\}^{2}$). 
Then, observe that the case depicted in Figure \ref{fig-BandM} corresponds with the matrix multiplication
\begin{equation}
\label{eq-ShiftExample}
\overbrace{ \begin{pmatrix}
E_{2}, E_{2}
\end{pmatrix} }^{ \TransBlock{M}{1}{1} }
\overbrace{ \begin{pmatrix}
E_{4}, E_{4}
\end{pmatrix} }^{ \TransBlock{M}{2}{2} \TransBlock{C}{3}{1}}
\overbrace{ 
\left(
\begin{array}{cc}
\begin{matrix}
\TransBlock{q}{2}{1} &  \\
 & \TransBlock{q}{2}{2}
\end{matrix}
& \text{\Large{\:0}} \\
\text{\Large{\:0}} 
& 
\begin{matrix}
\TransBlock{q}{2}{3} &  \\
 & \TransBlock{q}{2}{4}
\end{matrix}
\end{array}
\right)
}^{ \TransBlock{B}{2}{2}\left( \TransBlock{\mathcal{Q}}{2}{2} \right)}
\overbrace{
\begin{pmatrix}
\TransBlock{q}{1}{1} & \zero{2} \\
\zero{2} & \TransBlock{q}{1}{2}
\end{pmatrix}
}^{ \TransBlock{B}{1}{1}\left( \TransBlock{\mathcal{Q}}{1}{2} \right) }
= 
\sum_{i=0,1}
\left(
\begin{array}{cc}
 \begin{matrix}
q_{0, i, 0 \mid  0} & q_{1, i, 0 \mid 1} \\
q_{0, i, 1 \mid 0} & q_{1, i, 1 \mid 1} 
\end{matrix} 
\end{array}
\right),
\end{equation}
where $q_{k, i, j \mid j} = q_{k, i, j \mid i, j} q_{i, j \mid j}$.
The right-hand side of (\ref{eq-ShiftExample}) represents transition probability from $P(X_{t})$ to $P(X_{t+2})$ after the ``marginalization'' of $P(X_{t+1})$.
This is a special case of the so-called shift matrix $\mathcal{S}\left( \Trans{ \mathcal{Q} }{1}, \Trans{ \mathcal{Q} }{2} \right)$, 
which will be defined shortly.

The reader can confirm the properties of these operators above by finding the following identities.
For an arbitrary tuple of simplex vectors $\TransBlock{\mathcal{Q}}{k}{N}$ and $0 \le m, n \le k$, we have the identity
\[
 \TransBlock{M}{k}{m} = \TransBlock{C}{k}{m-n} \TransBlock{M}{k}{n} \TransBlock{C}{k+1}{n-m},
\]
\[
 \TransBlock{M}{k}{n} \TransBlock{M}{k+1}{m} = 
\begin{cases}
 \TransBlock{M}{k}{m} \TransBlock{M}{k+1}{n+1} \ \text{for}\ n \ge m \\
 \TransBlock{M}{k}{m-1} \TransBlock{M}{k+1}{n} \ \text{otherwise},
\end{cases}
\]
and
\[
 \TransBlock{B}{k}{m}\left( \TransBlock{\mathcal{Q}}{k}{N} \right) = 
 \TransBlock{C}{k+1}{m-n} \TransBlock{B}{k}{n}\left( \TransBlock{\mathcal{Q}}{k}{N} \right) \TransBlock{C}{k}{n-m}.
\]
For an arbitrary integer $m, m$, 
\[
 \TransBlock{C}{k}{ \mathrm{mod}_{k}\left( n m \right) } = \left( \TransBlock{C}{k}{n} \right)^{m}
\]

Now it is easy to understand that the $k^{\text{th}}$ order transition matrix $\Trans{Q}{k}$ (\ref{eq-TransitionMat}), 
which is the ``shift'' operator $P\left( X_{1}^{k} \right) \rightarrow P\left( X_{2}^{k+1} \right)$, can be written with the three matrices as follows.

\begin{proposition}[Transition matrix as the shift operator]
\label{thm-Diagonalization}
Denote an arbitrary transition matrix 
by $\Trans{Q}{k} \in \SimMat{N^{k}}$, and the corresponding tuple of vectors 
$\TransBlock{\mathcal{Q}}{k}{N} := \left( \TransBlock{q}{k}{N(i-1)+j} \right)_{
i \in \Encode{N}{k-1}, j \in \Encode{N}{1}}$. Then we have
\[
\Trans{Q}{k} = \TransBlock{M}{k}{k} \TransBlock{C}{k+1}{1} \TransBlock{B}{k}{k}\left( \TransBlock{\mathcal{Q}}{k}{N} \right)
 = \TransBlock{C}{k}{1} \TransBlock{M}{k}{k-1} \TransBlock{B}{k}{k}\left( \TransBlock{\mathcal{Q}}{k}{N} \right)
.\]
\end{proposition}
\begin{proof}
\begin{eqnarray}
\TransBlock{M}{k}{k-1} \TransBlock{B}{k}{k}\left( \TransBlock{\mathcal{Q}}{k}{N} \right)  &=& 
 \sum_{i \in \Encode{N}{k-1}, j \in \Encode{N}{1}} E_{N^{k-1},i} 
         \otimes e_{N, j}^{T}   \otimes \TransBlockVec{q}{k}{i}{j}
\nonumber
\\
\nonumber
&=& \sum_{i \in \Encode{N}{k-1}} E_{N^{k-1}, i} \otimes \TransBlock{Q}{k}{i} = C_{N, N^{k-1}}^{-1} \TransBlock{Q}{k}.
\end{eqnarray}
\end{proof}


\subsection{$k$-shift matrix \label{sec-shift}}
As we wish to calculate the marginal distribution of $k^{\text{th}}$ Markov process 
rather than the full stationary distribution, 
it is crucial to describe the property of the marginal stationary distribution.
Here, 
we will define the {\em $k$-shift matrix}, with which 
 we analyze the recursive property of an arbitrary marginal stationary distribution.


Let us introduce a $k$-shift operator, which is 
 an extension of Proposition \ref{thm-Diagonalization}, as follows.

\begin{definition}
\label{def-kshift}
For a series of tuples $\Trans{\mathcal{Q}}{1}, \ldots, \Trans{ \mathcal{Q} }{k}$, 
define a $k$-shifting transition matrix:
\[
 \mathcal{S}\left( \Trans{\mathcal{Q}}{1}, \ldots, \Trans{ \mathcal{Q} }{k} \right) := 
\TransBlock{M}{1}{1}
 \ldots 
\TransBlock{M}{k}{k}\TransBlock{C}{k+1}{1}
          \TransBlock{B}{k}{k}\left( \Trans{\mathcal{Q}}{k} \right) 
   \ldots \TransBlock{B}{1}{1}\left( \Trans{\mathcal{Q}}{1} \right).
\]
\end{definition}

We can easily see that the transition matrix is identical to the 1-shift matrix, 
$\Trans{Q}{k} = \mathcal{S}\left( \Trans{ \mathcal{Q} }{k} \right)$, by this definition.
Find an example of 2-shift matrix $\mathcal{S}\left( \Trans{ \mathcal{Q} }{1}, \Trans{ \mathcal{Q} }{2} \right)$ in (\ref{eq-ShiftExample}).
As seen in (\ref{eq-ShiftExample}), the $k$-shift matrix ``pushes'' the marginal probability vector of $P(X_{t})$ to that of $P(X_{t+k})$.
Thus, generally a $k$-shift matrix expresses transition probability of the marginalized states for a given Markov process. 

In what follows, we show an important property of the $k$-shift matrix, namely that 
it is written in a recursive form.
Let us consider the same $\mathcal{S}\left( \Trans{ \mathcal{Q} }{1}, \Trans{ \mathcal{Q} }{2} \right)$ as in (\ref{eq-ShiftExample}).
It has another compact form of expression
\[
\mathcal{S}\left( \Trans{ \mathcal{Q} }{1}, \Trans{ \mathcal{Q} }{2} \right)
= \left( \TransBlock{\overline{Q}}{2}{1} \TransBlock{q}{1}{1}, \TransBlock{\overline{Q}}{2}{2} \TransBlock{q}{1}{2} \right),
\]
with $\TransBlock{\overline{Q}}{2}{1} := \left( \TransBlock{q}{2}{1}, \TransBlock{q}{2}{2} \right)$
and 
$\TransBlock{\overline{Q}}{2}{2} := \left( \TransBlock{q}{2}{3}, \TransBlock{q}{2}{4} \right)$.

Similarly, with an additional tuple $\Trans{ \mathcal{Q} }{3} := \left( \TransBlock{q}{3}{1}, \ldots, \TransBlock{q}{3}{2^{3}} \right)$, observe 
\[
\mathcal{S}\left( \Trans{ \mathcal{Q} }{1}, \Trans{ \mathcal{Q} }{2}, \Trans{ \mathcal{Q} }{3} \right)
= \left( \TransBlock{\overline{Q}}{2}{1} \TransBlock{q}{1}{1}, \TransBlock{\overline{Q}}{2}{2} \TransBlock{q}{1}{2} \right),
\]
with $\TransBlock{\overline{Q}}{2}{i} := \left( \TransBlock{\overline{Q}}{3}{2(i-1)+1}\TransBlock{q}{2}{2(i-1)+1}, \TransBlock{\overline{Q}}{3}{2(i-1)+1}\TransBlock{q}{2}{2(i-1)+1} \right)$ for $i=1,2$
and 
$\TransBlock{\overline{Q}}{3}{i} := \left( \TransBlock{q}{3}{2(i-1)+1}, \TransBlock{q}{3}{2(i-1)+1} \right)$ for $i=1, \ldots, 4$.

In general, as expected from these examples above, 
a $k$-shift matrix can be written in a recursive form.
With $ \TransBlock{\overline{Q}}{m+1}{i} := E_{2}$ for any $1 \le i \le 2^{m}$, 
$\mathcal{S}\left( \Trans{ \mathcal{Q} }{1}, \ldots, \Trans{ \mathcal{Q} }{m} \right) = \TransBlock{ \overline{Q} }{1}{1}$ is defined recursive as follows: 
\begin{equation}
\nonumber
 \begin{matrix}
 \underbrace{
 \left( 
      \TransBlock{\overline{Q}}{m+1}{1} \TransBlock{q}{m}{1}, 
      \TransBlock{\overline{Q}}{m+1}{2} \TransBlock{q}{m}{2}
  \right) 
 }_{=\TransBlock{\overline{Q}}{m}{1}},
 & \ldots, & 
  \underbrace{
 \left( 
      \TransBlock{\overline{Q}}{m+1}{2^{m}-1} \TransBlock{q}{m}{2^{m}-1}, 
      \TransBlock{\overline{Q}}{m+1}{2^{m}} \TransBlock{q}{m}{2^{m}}
  \right) 
 }_{=\TransBlock{\overline{Q}}{m}{2^{m-1}}}
\\
 \underbrace{
 \left( 
      \TransBlock{\overline{Q}}{m}{1} \TransBlock{q}{m-1}{1}, 
      \TransBlock{\overline{Q}}{m}{2} \TransBlock{q}{m-1}{2}
  \right) 
 }_{=\TransBlock{\overline{Q}}{m-1}{1}},
 & \ldots, & 
  \underbrace{
 \left( 
      \TransBlock{\overline{Q}}{m}{2^{m-1}-1} \TransBlock{q}{m-1}{2^{m-1}-1}, 
      \TransBlock{\overline{Q}}{m}{2^{m-1}} \TransBlock{q}{m-1}{2^{m-1}}
  \right) 
 }_{=\TransBlock{\overline{Q}}{m-1}{2^{m-2}}}
\\
& \vdots &
\\
 \left( 
      \TransBlock{\overline{Q}}{2}{1} \TransBlock{q}{1}{1}, 
      \TransBlock{\overline{Q}}{2}{2} \TransBlock{q}{1}{2}
  \right) 
={\TransBlock{\overline{Q}}{1}{1}} && 
 \end{matrix}
\end{equation}

For a general $N$, the following lemma asserts that this recursive property holds.
\begin{lemma}[Recursive property of the $k$-shift matrix]
\label{prop-RecursiveForm}
Given tuples of vectors 
$\TransBlock{\mathcal{Q}}{m}{N} := \left( \TransBlock{q}{m}{1}, \ldots, \TransBlock{q}{m}{N^{m}} \right) \in \SimMats{N}{N^{m}}$
for $m = 0, 1, \ldots, k$,
we can write the corresponding $k$-shifting transition matrix in a recursive form:
\[
 \mathcal{S}\left(   \TransBlock{\mathcal{Q}}{m}{N}, \ldots, \TransBlock{\mathcal{Q}}{k}{N} \right) 
= \TransBlock{C}{m}{1}\sum_{i=1}^{N^{m-1}} E_{N^{m-1},i} \otimes \TransBlock{\overline{Q}}{m}{i}
\]
where
for $1 \le m \le k$ and $1 \le i \le N^{m-1}$
\[
 \TransBlock{\overline{Q}}{m}{i}
:=
\left( \TransBlock{\overline{Q}}{m+1}{N(i-1)+1} \TransBlockVec{q}{m}{i}{1}, 
 \TransBlock{\overline{Q}}{m+1}{N(i-1)+2} \TransBlockVec{q}{m}{i}{2}, 
 \ldots, 
 \TransBlock{\overline{Q}}{m+1}{N(i-1)+N} \TransBlockVec{q}{m}{i}{N}
  \right),
\]
and 
$\TransBlock{\overline{Q}}{k+1}{i} := E_{N}$ for $1 \le i \le N^{k}$.
\end{lemma}
\begin{proof}
Observe the recurrent relationship between
\[
 \mathcal{S}\left( \TransBlock{\mathcal{Q}}{k}{N} \right)
 = \Trans{Q}{k} = 
   \sum_{i=1}^{N^{k-1}} e_{N^{k-1},i}^{T} \otimes \TransBlock{Q}{k}{i} \otimes e_{N^{k-1},i}
\]
and
\begin{eqnarray}
\mathcal{S}\left( \TransBlock{\mathcal{Q}}{k-1}{N}, \TransBlock{\mathcal{Q}}{k}{N} \right)
 & = & \TransBlock{C}{k-1}{1} \TransBlock{M}{k-1}{k-1} 
  \left(\TransBlock{C}{k}{1}\right)^{-1} \Trans{Q}{k}
  \TransBlock{B}{k-1}{k-1}\left(  \TransBlock{\mathcal{Q}}{k}{N} \right)
\nonumber 
\\
&=& \sum_{i=1}^{N^{k-2}} e_{N^{k-2},i}^{T} \otimes 
  \TransBlock{ \overline{Q} }{k-1}{i}
  \otimes e_{N^{k-2},i}
\nonumber
\end{eqnarray}
where 
$\TransBlock{ \overline{Q} }{k-1}{i} = \sum_{j=1}^{N} e_{N,j}^{T} \otimes \TransBlock{ {Q} }{k}{i}\TransBlockVec{q}{k-1}{i}{j}$.
For $1 \le m < k$, 
find the recursive relationship between
$\mathcal{S}\left( \TransBlock{\mathcal{Q}}{m-1}{N}, \ldots, \TransBlock{\mathcal{Q}}{k}{N} \right)$
and $\mathcal{S}\left( \TransBlock{\mathcal{Q}}{m}{N}, \ldots, \TransBlock{\mathcal{Q}}{k}{N} \right)$
by inductively writing 
$\TransBlock{ \overline{Q} }{m-1}{i} := \sum_{j=1}^{N} e_{N,j}^{T} \otimes \TransBlock{ Q }{m}{i}
\TransBlockVec{q}{m-1}{i}{j}$.
\end{proof}

\subsection{Recursive property of marginal stationary distribution}
Using the $k$-shift matrix, here we will observe the recursive property of the marginal stationary distribution.
Before that, however, let us briefly remark that an arbitrary probability vector is uniquely expressed with branching matrices as follows.

\begin{proposition}
\label{prop-1to1}
 For every $\theta \in \Simplex{N^{k}}$ without any zero element, there is a unique tuple of vectors  
$\TransBlock{{\Theta}}{m}{N} := \left( \TransBlock{\theta}{m}{1}, \ldots, \TransBlock{\theta}{m}{N^m} \right)_{
i \in \Encode{N}{m}}$ for $m = 0, 1, \ldots, k-1$, which holds
\[
 \theta = \TransBlock{C}{k}{1}
          \TransBlock{B}{k-1}{k-1}\left( \TransBlock{{\Theta}}{k-1}{N} \right) 
          \TransBlock{B}{k-2}{k-2}\left( \TransBlock{{\Theta}}{k-2}{N} \right) 
   \ldots \TransBlock{B}{0}{0}\left( \TransBlock{{\Theta}}{0}{N} \right).
\]
\end{proposition}
\begin{proof}
As $ \TransBlock{M}{k}{k}\TransBlock{B}{k}{k}\left( \TransBlock{\Theta}{k}{N} \right) = E_{N^{k-1}}$ for an arbitrary $k$,
define 
$
\TransBlock{\theta}{0}{1} := \TransBlock{M}{1}{1} \ldots \TransBlock{M}{k-1}{k-1} \theta
$.
For $i > 0$ and $\TransBlock{\Theta}{i-1}{N}:= \left( \TransBlock{\theta}{i-1}{1}, \dots, \TransBlock{\theta}{i-1}{N^{i-1}} \right)$, 
define $\TransBlock{\theta}{i}{1}, \dots, \TransBlock{\theta}{i}{N^{i}} $ 
by the root of 
\begin{equation}
\label{eq-Marginal}
\TransBlock{C}{i+1}{1}
\overbrace{
\begin{pmatrix}
\TransBlock{\theta}{i}{1} & & \text{\Large{\:0}} 
\\
 & \ddots & 
\\
\text{\Large{\:0}}  & & \TransBlock{\theta}{i}{N^{i}}
\end{pmatrix}
}^{ \TransBlock{B}{i}{i}\left( \TransBlock{\Theta}{i}{N} \right) }
 \TransBlock{B}{i-1}{i-1}\left( \TransBlock{\Theta}{i-1}{N} \right) \ldots \TransBlock{B}{0}{0}\left( \TransBlock{\Theta}{0}{N} \right) 
= \TransBlock{M}{i+1}{i+1} \ldots \TransBlock{M}{k-1}{k-1} \theta.
\end{equation}
For $\TransBlock{\Theta}{i-1}{N}$ without any zero element, which is necessary for $\theta$ without any zero element, this gives a unique $\TransBlock{\Theta}{i}{N} := \left( \TransBlock{\theta}{i}{1}, \dots, \TransBlock{\theta}{i}{N^{i}} \right)$.
Repeat this construction up to $i = k-1$, and the right hand side of (\ref{eq-Marginal}) is $\theta$.
\end{proof}

We are now ready to state the lemma on the recursive property of the marginal stationary distribution.

\begin{lemma}
\label{prop-marginal}
For a transition matrix $\Trans{Q}{k} \in \SimMat{N^{k}}$, let its unique stationary vector be $\theta = \Trans{Q}{k} \theta \in \Simplex{N^{k}}$,
and suppose it is an irreduceable Markov process or equivalently $\theta$ has no zero element.
From Proposition \ref{prop-1to1}, 
write 
$$ \theta = \TransBlock{C}{k}{1} \Trans{B}{k-1}\left( \TransBlock{{\Theta}}{k-1}{N} \right) 
   \ldots \Trans{B}{0}\left( \TransBlock{{\Theta}}{0}{N} \right)
\in \Simplex{N^{k}}$$ 
with
$\TransBlock{{\Theta}}{m}{N} := \left( \TransBlock{\theta}{m}{1}, \ldots, \TransBlock{\theta}{m}{N^m} \right)$ for $m = 0, 1, \ldots, k$.
Denote the $m^{\text{th}}$ order marginalized stationary vector by 
$$\Trans{\theta}{m} := \TransBlock{C}{m}{1} \TransBlock{B}{m-1}{m-1}\left( \TransBlock{\Theta}{m-1}{N} \right) \ldots 
\TransBlock{B}{0}{0}\left( \TransBlock{\Theta}{0}{N} \right) \in \Simplex{N^{m}}.$$
Then, we have for any $m \le i \le k$
\[
 \Trans{\theta}{m} = \mathcal{S}\left( \Trans{\Theta}{m}, \ldots, \Trans{\Theta}{i} \right) \Trans{\theta}{m} = 
\mathcal{S}\left( \Trans{\Theta}{m}, \ldots, \Trans{\Theta}{k}, \Trans{Q}{k} \right) \Trans{\theta}{m}.
\]
\end{lemma}
\begin{proof}
For $m \le i \le k$, 
multiply $\TransBlock{M}{m}{m} \ldots \TransBlock{M}{i}{i}$
from the left to (\ref{eq-Marginal}) in Proposition \ref{prop-1to1}, and we have
\begin{equation}
\nonumber
 \mathcal{S}\left(   \TransBlock{ \Theta }{m}{N}, \ldots, \TransBlock{ \Theta }{i}{N} \right) 
\Trans{\theta}{m}
= 
\TransBlock{M}{m}{m}
\ldots
\TransBlock{M}{k-1}{k-1} \theta
= \Trans{\theta}{m}.
\end{equation}
The left most one above is due to Definition \ref{def-kshift}, and the rightmost one above is due to (\ref{eq-Marginal}).
Applying Proposition \ref{thm-Diagonalization}, $\Trans{Q}{k} = \TransBlock{M}{k}{k}\TransBlock{C}{k+1}{1}\TransBlock{B}{k}{k}\left( \TransBlock{ \mathcal{Q} }{k}{N} \right)$ to $\theta = \Trans{Q}{k} \theta$, 
 we have
$$\Trans{\theta}{m} =  \mathcal{S}\left(   \Trans{ \Theta }{m}, \ldots, \Trans{ \Theta }{k}, \Trans{Q}{k} \right)  \Trans{\theta}{m}.$$
\end{proof}

\subsection{$k$-shift matrix and marginal distribution}

Lemma \ref{prop-marginal} shows the essential similarity between stationary vector of $k$-shift matrix and marginalized stationary vector of the corresponding $k^{\text{th}}$ order Markov process.
Let us denote the $m^{\text{th}}$ marginal vector of the $k^{\text{th}}$ order Markov process by $\TransBlock{\theta}{m}{k}$, and it holds $\TransBlock{\theta}{k}{k-1} = \Trans{Q}{k}\TransBlock{\theta}{k}{k-1}$ by definition.
The marginal stationary vector of the $k^{\text{th}}$ order Markov process then holds
\[
 \TransBlock{\theta}{0}{k} =  \mathcal{S}\left( \TransBlock{\Theta}{1}{k}, \ldots, \TransBlock{\Theta}{k}{k}, \Trans{ \mathcal{Q} }{k} \right)\TransBlock{\theta}{0}{k},
\]
while the stationary vector $\omega_{k}$ of the $k$-shift matrix holds
\[
 \omega_{k} =  \mathcal{S}\left( \Trans{\mathcal{Q}}{1}, \ldots, \Trans{ \mathcal{Q} }{k} \right)\omega_{k}.
\]
The following theorem asserts a relationship between the two types of stationary vectors.

\begin{theorem}
\label{thm-RecursiveMarginal}
Consider a series of tuples of vectors $\Trans{ \mathcal{Q} }{1}, \Trans{ \mathcal{Q} }{2}, \ldots$
and let $\Trans{ Q }{1}, \Trans{ Q }{2}, \ldots$ and  $\Trans{ \theta }{1}, \Trans{ \theta }{2}, \ldots$ 
be the corresponding series of transition matrices and stationary vectors, respectively. Suppose each forms an irreduceable Markov process 
for every $k$, and 
the series of the $m^{\text{th}}$ marginal of the $k^{\text{th}}$ order stationary vectors $\TransBlock{\theta}{m}{k}$ 
has the limit for any $0 \le m < k$ and an arbitrary integer $l \ge 0$
\begin{equation}
\label{eq-ThetaCondition}
 \lim_{k \to \infty} \left\| \TransBlock{\theta}{m}{k} - \TransBlock{\theta}{m}{k+l} \right\| = 0.
\end{equation}
Suppose that each $m^{\text{th}}$ order marginal stationary vector of the $k$-shift matrix for $0\le m < k $, 
$$
\TransBlock{\omega}{m}{k} = 
\mathcal{S}\left( \Trans{ \mathcal{Q} }{m}, \ldots, \Trans{ \mathcal{Q} }{k} \right) \TransBlock{\omega}{m}{k},
$$
 has a unique root and forms a converging series for $l \ge 0$
\begin{equation}
\label{eq-OmegaCondition}
 \lim_{k \to \infty} \left\| \TransBlock{\omega}{m}{k} - \TransBlock{\omega}{m}{k+l} \right\| = 0.
\end{equation}
We then have 
\[
 \lim_{k \to \infty}\left\| \TransBlock{\theta}{m}{k} - \TransBlock{\omega}{m}{k} \right\| = 0.
\]
\end{theorem}
\begin{proof}
Condition (\ref{eq-ThetaCondition}) implies that there is some $k_{0}$ for an arbitrary $\epsilon > 0$ which lets $\left\| \TransBlock{\theta}{m}{k} - \TransBlock{\theta}{m}{k+1} \right\| < \epsilon$ for $k \ge k_{0}$.
According to Lemma \ref{prop-marginal}, this implies 
$$
\left\| \left( 
\mathcal{S}\left( \TransBlock{\Theta}{m}{k}, \ldots, \TransBlock{\Theta}{k}{k}, \Trans{ \mathcal{Q} }{k} \right) 
- 
\mathcal{S}\left( \TransBlock{\Theta}{m}{k+1}, \ldots, \TransBlock{\Theta}{k}{k+1}, \TransBlock{ \Theta }{k+1}{k+1} \right) 
\right)\TransBlock{\theta}{m}{k+1} \right\| < \epsilon,
$$
and by the Cauchy-Schwartz inequality we have 
$$
\left\| \left( 
\mathcal{S}\left( \TransBlock{\Theta}{m}{k}, \ldots, \TransBlock{\Theta}{k}{k}, \Trans{ \mathcal{Q} }{k}, \Trans{ \mathcal{Q} }{k+1} \right) 
- 
\mathcal{S}\left( \TransBlock{\Theta}{m}{k}, \ldots, \TransBlock{\Theta}{k}{k+1}, \TransBlock{ \Theta }{k+1}{k+1}, \Trans{ \mathcal{Q} }{k+1} \right) 
\right)\TransBlock{\theta}{m}{k+1} \right\| < \epsilon.
$$
Apply this inequality sequentially to $k+2, \ldots, k+l$, and we have 
$$
\left\| \left( 
E_{N^{m+1}}
- 
\mathcal{S}\left( \TransBlock{\Theta}{m}{k}, \ldots, \TransBlock{\Theta}{k}{k}, \Trans{ \mathcal{Q} }{k}, \ldots, \Trans{ \mathcal{Q} }{k+l} \right) 
\right)\TransBlock{\theta}{m}{k+l} \right\| < \epsilon.
$$
For $\TransBlock{ \omega }{ m }{ k }$, which holds the condition (\ref{eq-OmegaCondition}), 
we have 
\[
 \lim_{k \to \infty}\left\| \TransBlock{\theta}{m}{k} - \TransBlock{\omega}{m}{k} \right\| = 0,
\]
due to the irreduceability of the series of Markov processes, which makes the convergence unique. 
\end{proof}

\section{Recursive Markov process \label{sec-Recursive}}
Theorem \ref{thm-RecursiveMarginal} implies that 
the marginal stationary vector $\TransBlock{\theta}{k}{m}$
is closely approximated by stationary vector $\TransBlock{\omega}{k}{m}$ of the $k$-shift matrix
$\mathcal{S}\left( \TransBlock{\mathcal{Q}}{m}{N}, \ldots,  \TransBlock{\mathcal{Q}}{k}{N} \right)$ under the limit $k \to \infty$.
This theorem motivates us to consider a special class of Markov processes, which can be constructed in a recursive manner, as follows.


\begin{definition}[Recursive Markov process]
\label{def-Recursive}
 We call a $k^{\text{th}}$ order Markov process with the transition matrix $\Trans{Q}{k}$ {\em recursive},
if each element of the block matrix $\TransBlock{Q}{m}{i}$ 
is a function of the elements of $\TransBlock{q}{m-1}{i} \in \Simplex{N}$
for 
$ 1 < m \le k$, $i \in \Encode{N}{m}$.
\end{definition}
This definition of the recursive Markov process is motivated by the fact that 
we can analyze the convergence of such a series of transition matrices in a closed form.
The following corollary states 
that this class is characterized by a closed-form equation of the marginal stationary distribution.
\begin{corollary}
\label{cor-recursive}
Suppose there is 
a map $f: \Simplex{N} \mapsto \SimMat{N}$, with which an infinite order recursive Markov process satisfies
$\TransBlock{Q}{m+1}{i} = f\left( \TransBlock{q}{m}{i} \right)$ for $m = 1, 2, \ldots $ and $i \in \Encode{N}{m}$. 
Denote the fixed point $ \omega \in \Simplex{N}$ for the linear transformation $f( \omega )$, 
which satisfies $ \omega = f( \omega ) \omega $.
Then, the marginal stationary vector of the $k^{\text{th}}$ order stationary vector
 $\Trans{\theta}{k}$ under the limit $k \to \infty$ corresponds with $\omega$ as follows:
\[
 \omega = \lim_{k\to\infty} \TransBlock{M}{1}{1} \ldots \TransBlock{M}{k-1}{k-1}\Trans{\theta}{k} \in \Simplex{N},
\]
if the limit shift matrix
$\overline{Q} := \lim_{k\to\infty}\mathcal{S}\left( \TransBlock{\mathcal{Q}}{1}{N}, \ldots, \TransBlock{\mathcal{Q}}{k}{N} \right)$
of this recursive Markov process is irreducible.
\end{corollary}
\begin{proof}
 Denote 
$
\overline{Q} = \left( q_{1}, \ldots, q_{N} \right) \in \SimMat{N}
$.
According to Theorem \ref{thm-RecursiveMarginal}, the marginal stationary distribution 
$\omega$ holds $\omega = \overline{Q} \omega$, and 
the recursive Markov process holds $q_{i} = \overline{Q} q_{i} = f(q_{i}) q_{i}$ for $i \in \Encode{N}{1}$.
As $\overline{Q}$ is irreducible, $\omega = q_{i}$ and $\omega = f(\omega) \omega$.
\end{proof}


\section{Numerical case study \label{sec-Application}}

To see an application of the current mathematical analysis, 
we calculated the stationary vector for the $m$-shift matrix 
for a specific case of the iterated prisoner's dilemma with learning (Definition \ref{def-LearningGames}).
Here we present an analysis of the iterated prisoners' dilemma with rewards $R_{00} = 1, R_{01} = -2, R_{10} = 2$, $R_{11} = 0$ (Definition \ref{def-LearningGames}). The players are reinforcement learners with the identical sensitivity parameters $\beta_{1} = \beta_{2} = \beta = 1/2$ or $1$ , and identical memory retention parameters for two players $\alpha_{1} = \alpha_{2} = \alpha \in [0, 1]$ (in Equation (\ref{eq-ChoiceIndividual})). 

\subsection{Finite $m$-shift matrix}

Using Proposition \ref{prop-RecursiveForm} for a finite $m$, we calculated 
the $m$-shift transition matrix and its stationary vector, and obtained the marginal stationary probabilities 
$\theta = ( P(\text{CC}), P(\text{CD}), P(\text{DC}), P(\text{DD}) )^{T}$
of mutual cooperation ($\text{CC}$), mutual defection ($\text{DD}$), and one-side defection ($\text{CD}$ and $\text{DC}$, as $P(CD) = P(DC)$ with equality due to the symmetry between two players).
Letting the parameter $m$ be large, we expect that
the stationary vector calculated by the eigenvector of the $m$-shift matrix would be close to the marginal stationary vector calculated for the $k^{\text{th}}$ order Markov process with a sufficiently large $k$. 
Due to the limitations of our computational resources, we calculated it up to $m=12$. The calculation of these stationary probabilities takes on the order $4^{m+1} \approx 10^{7.82}$ steps for $m =12$.
For comparison, we compute the approximation for the exchangeable special case in which $\alpha = 1$ and $k=1000$ (See also Section \ref{sec-Exchangeable}).

\begin{figure}
    \begin{center}
    \includegraphics[width=1\linewidth, clip]{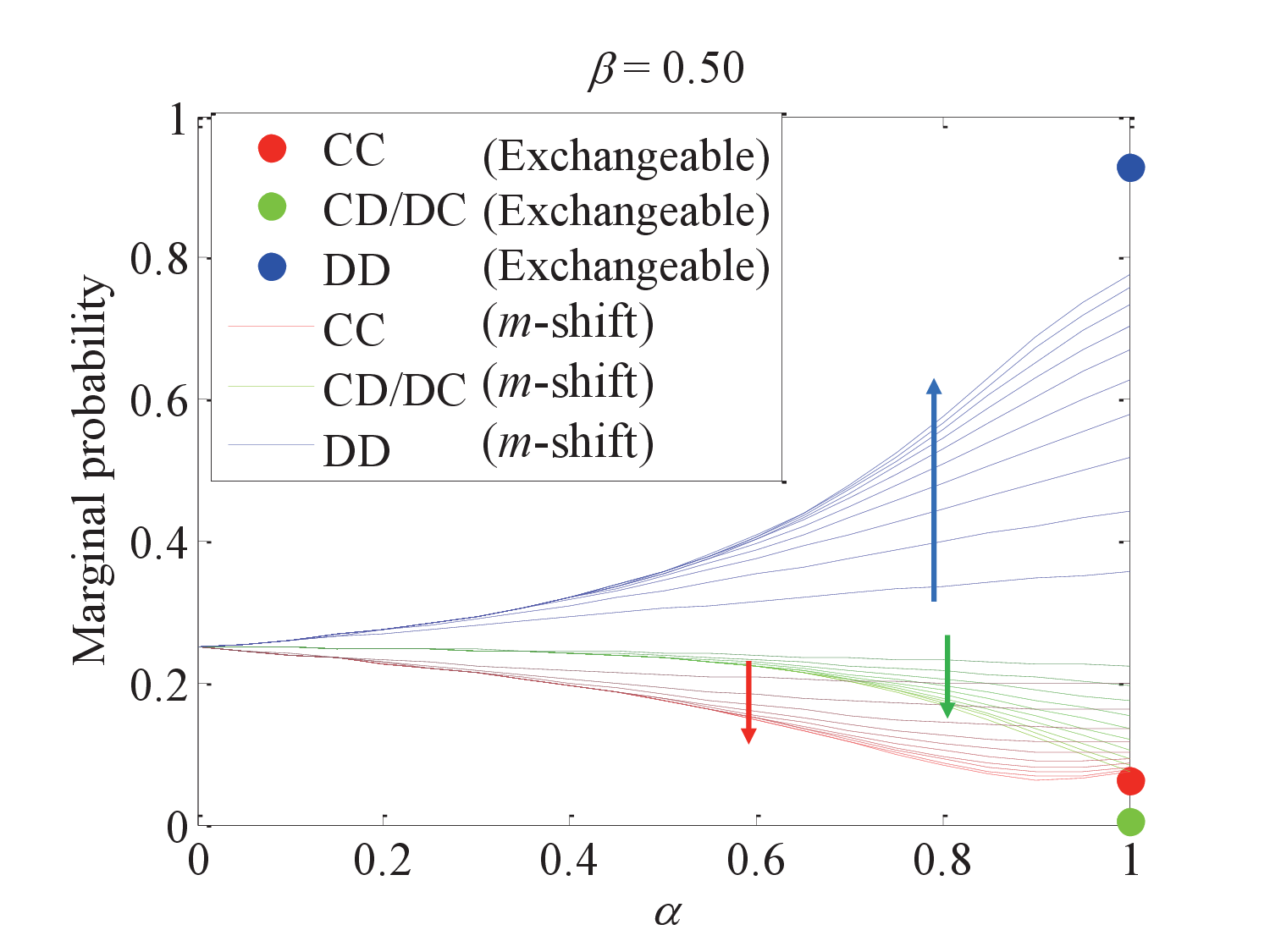}
    \caption{\label{fig-IPDMarginal} The $m$-shift stationary distributions (curves) and marginal distributions for $\alpha = 1$ and $k=1000$ (circles) as functions of the memory retention parameter $\alpha$.}
    \end{center}
  \end{figure}

Figure \ref{fig-IPDMarginal} shows the marginal stationary probabilities with $\beta = 1/2$ estimated by the $m$-shift stationary probabilities as functions of the memory retention parameter $\alpha$.
The multiple lines of the same color show the marginal probabilities for different values of the shift $m$ for a fixed value of $\alpha$.
The arrows indicate the directions in which these groups of lines change from $m = 1, 2, \hdots, 12$.

Observe that the probability for mutual defection $P(DD)$ increases with the memory retention parameter $\alpha$. This result is qualitatively consistent with the outcome of the classical prisoners' dilemma with two rational players.
Our analysis illustrates the counterpart of the classical Nash equilibrium in the iterated version of the game with probabilistic reasoners capable of remembering all the previous outcomes. Interpreting the memory retention parameter $\alpha$ as the degree of rationality of the agents, this indicates that, as players become more rational, they are more likely to mutually defect.

For the exchangable case $\alpha = 1$, the special computational procedure described in Section \ref{sec-Exchangeable} can be used even for large values of $k$. We used this procedure to perform an analysis of the exchangable case with $k = 1000$. In our analysis the estimated marginal distributions appear to converge. The results of this analysis are indicated by the filled circles in Figure \ref{fig-IPDMarginal}.

\begin{figure}
    \begin{center}
    \includegraphics[width=1\linewidth, clip]{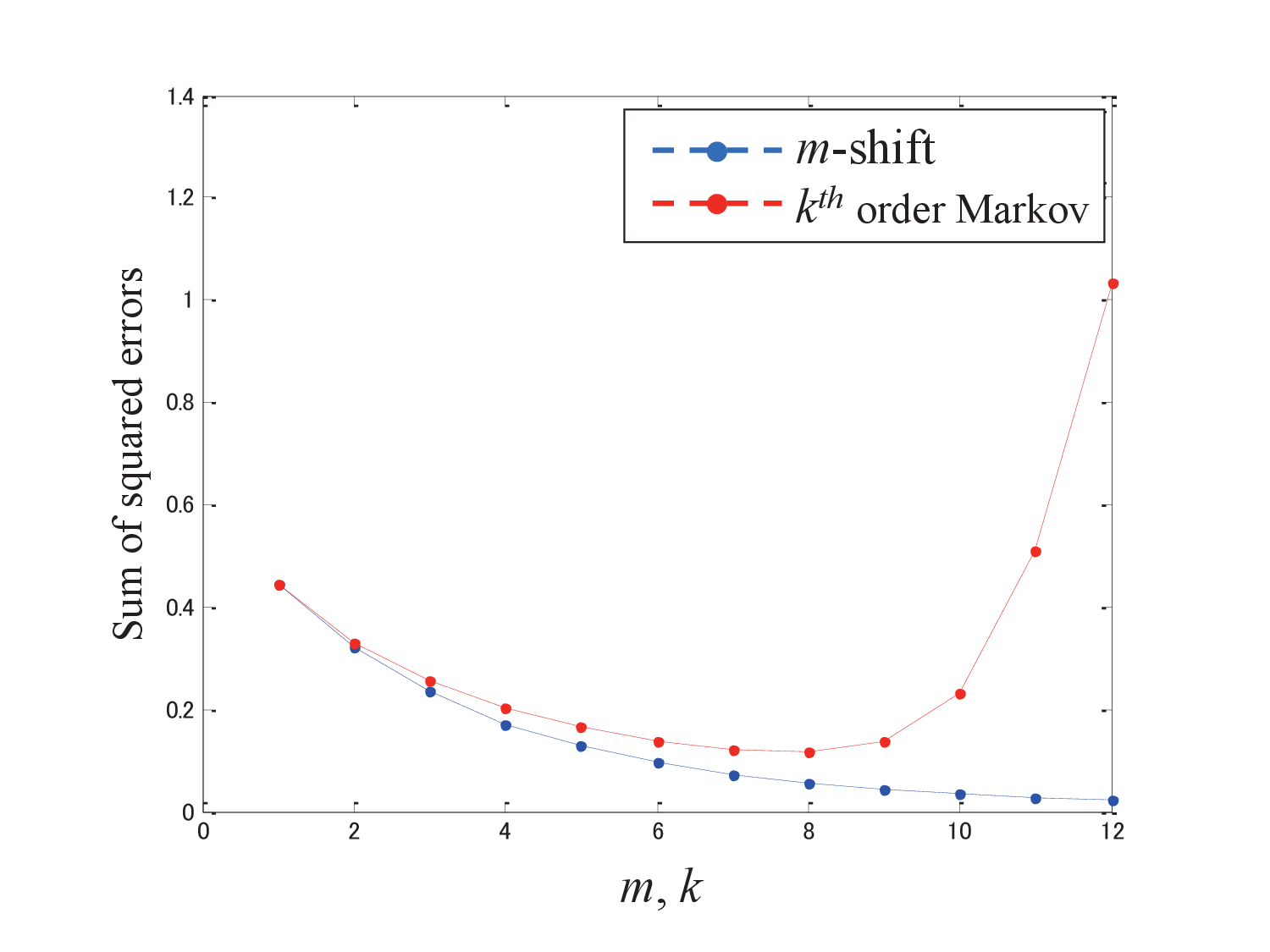}
    \caption{\label{fig-SSE}
The sum of squared errors of 
the $m$-shift stationary distributions (blue) and the $k^{\text{th}}$ order Markov process as functions of $m$ or $k$, respectively, by taking the corresponding marginal distribution for $\alpha = 1$ and $k=1000$ as normative values.}
    \end{center}
  \end{figure}

Treating the estimates obtained using the special property of $\alpha = 1$ as the true stationary distributions, we analyzed the sum of squared errors (SSE) in the estimated $m$-shift stationary distributions with $\alpha = 1$. The blue line in Figure \ref{fig-SSE} shows these SSEs as a function of $m$.
Consistent with Theorem \ref{thm-RecursiveMarginal}, the SSE is a decreasing function of $m$.
This error analysis shown in Figure \ref{fig-SSE} numerically endorses Theorem \ref{thm-RecursiveMarginal} by showing that the $m$-shift stationary distribution approaches a marginal stationary distribution under the limit $m \rightarrow \infty$.

In theory, as $k \rightarrow \infty$, the difference between the $k^{\text{th}}$-order stationary distributions and the corresponding $m$-shifted marginal stationary distributions could vanish in the limit. Such convergence, however, is not obvious.
The red line in Figure \ref{fig-SSE} shows SSE of the marginal distributions calculated by the $k^{\text{th}}$ order Markov process.
This result suggests that the errors of the $k^{\text{th}}$ order Markov process coverges more slowely (higher errors at each $k=m$) than the corresponding $m$-shift stationary distribution. 


\subsection{$\infty$-shift matrix}
Next, we additionally computed the stationary vector for the $\infty$-shift matrix by using Corollary \ref{cor-recursive}.
Specifically, 
we numerically solve the following non-linear equation for the vector $x = (x_{1}, x_{2}, x_{3}, x_{4})^{T}$:
\[
 x = M( x ) x,
\]
where
\[
 M(x) 
:= \left(
    \begin{array}{cccc}
      (a_{1}^{+}b_{1}^{+}x_{1})^{\alpha}/ c_{1} & 
      (a_{2}^{+}b_{2}^{-}x_{1})^{\alpha}/ c_{2} & 
      (a_{3}^{-}b_{3}^{+}x_{1})^{\alpha}/ c_{3} & 
      (a_{4}^{-}b_{4}^{-}x_{1})^{\alpha}/ c_{4} \\
      (a_{1}^{+}b_{1}^{-}x_{2})^{\alpha}/ c_{1} & 
      (a_{2}^{+}b_{2}^{+}x_{2})^{\alpha}/ c_{2} & 
      (a_{3}^{-}b_{3}^{-}x_{2})^{\alpha}/ c_{3} & 
      (a_{4}^{-}b_{4}^{+}x_{2})^{\alpha}/ c_{4} \\
      (a_{1}^{-}b_{1}^{+}x_{3})^{\alpha}/ c_{1} & 
      (a_{2}^{-}b_{2}^{-}x_{3})^{\alpha}/ c_{2} & 
      (a_{3}^{+}b_{3}^{+}x_{3})^{\alpha}/ c_{3} & 
      (a_{4}^{+}b_{4}^{-}x_{3})^{\alpha}/ c_{4} \\
      (a_{1}^{-}b_{1}^{-}x_{4})^{\alpha}/ c_{1} & 
      (a_{2}^{-}b_{2}^{+}x_{4})^{\alpha}/ c_{2} & 
      (a_{3}^{+}b_{3}^{-}x_{4})^{\alpha}/ c_{3} & 
      (a_{4}^{+}b_{4}^{+}x_{4})^{\alpha}/ c_{4} \\
    \end{array}
  \right)
\]
and
$c_{i}:= (a_{i}^{+}b_{i}^{+}x_{1})^{\alpha} + (a_{i}^{+}b_{i}^{-}x_{2})^{\alpha} + (a_{i}^{-}b_{i}^{+}x_{3})^{\alpha} + (a_{i}^{-}b_{i}^{-}x_{4})^{\alpha}$ under the constraint of $\sum_{i=1}^{4}x_{i} = 1$ and $x_{i} \ge 0$.
The set of parameters reflect the learning update of the choice probability (Equation (\ref{eq-ChoiceIndividual})) as the function of the payoffs: 
$a_{i}^{-}=b_{i}^{-} = 1$ for any $i$ and 
\begin{eqnarray}
a_{1}^{+} = \exp( R_{00} \beta ), 
a_{2}^{+} = \exp( R_{01} \beta ),
a_{3}^{+} = \exp( R_{10} \beta ),
a_{4}^{+} = \exp( R_{11} \beta ) \nonumber \\
b_{1}^{+} = \exp( R_{00} \beta ),
b_{2}^{+} = \exp( R_{10} \beta ),
b_{3}^{+} = \exp( R_{01} \beta ),
b_{4}^{+} = \exp( R_{11} \beta ) \nonumber .
\end{eqnarray}



\begin{figure}
    \begin{center}
    \includegraphics[width=1\linewidth, clip]{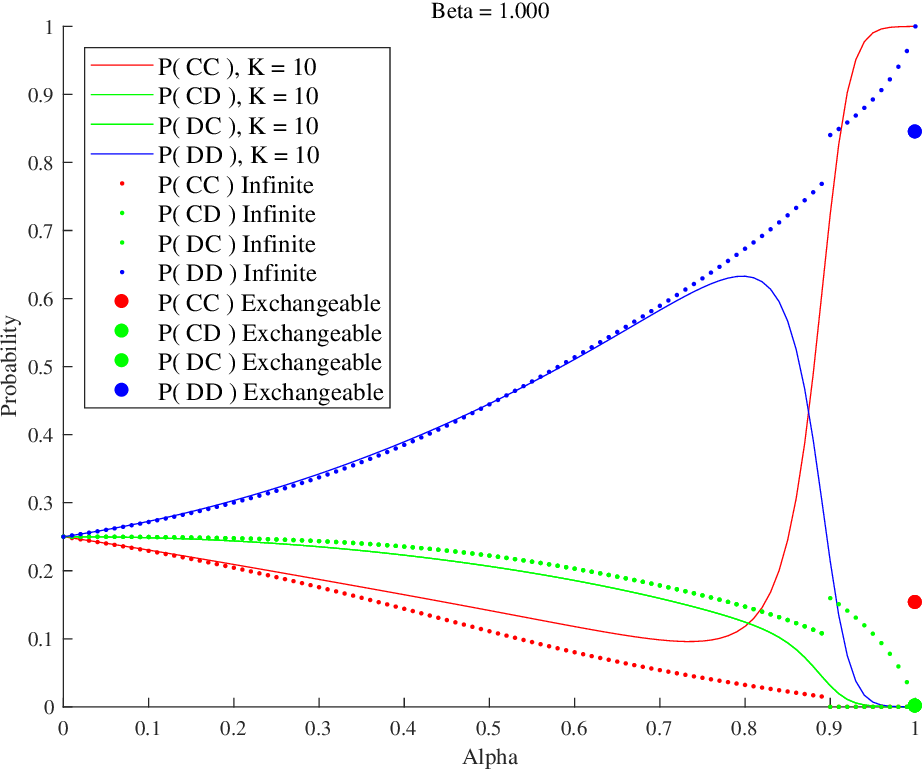}
    \caption{\label{fig-FiniteInfinite}
The marginal probability of the four pairwise states $( \text{CC}, \text{CD}, \text{DC}, \text{DD} )$ calculated by the 
$k^{\text{th}}$ order Markov process (lines), $\infty$-shift matrix (dots), and the exchageable special case for $\alpha = 1$ and $k=2000$ (large filled circles).}
    \end{center}
  \end{figure}

Figure \ref{fig-FiniteInfinite} shows the marginal stationary vector calculated by the $10^{\text{th}}$ order Markov process (lines), $\infty$-shift matrix (dots), and the exchangeable case for $\alpha = 1$.
The marginal stationary vector of $\infty$-shift matrix is expected to close to a finite order Markov process for a sufficiently small $\alpha$, and we observe the close match between those of the finite order Markov and $\infty$-shift matrix 
up to $\alpha \le 0.7$ in Figure \ref{fig-FiniteInfinite}. 
In addition, 
the marginal stationary vector of the $\infty$-shift matrix is expected to be close to that of the approximated one calculated by the procedure for the exchangeable case at $\alpha = 1$. As expected, we find the exchangeable one $k=1000$ is closer to that of the $\infty$-shift matrix than that of $10^{\text{th}}$ order Markov process. 
We find the gap of the probabilities of the stationary vector of the $\infty$-shift matrix at $\alpha = 0.9$, which is perhaps due to some numerical error in the calculation.

\subsection{Discussion}
In summary, this numerical study endorses the theoretically expected properties of the stationary vector of the $m$-shift matrix for both finite and infinite $m$.
Beyond the purpose of numerical validation, 
the result in Figure \ref{fig-FiniteInfinite} suggests that (1) a pair of players learning from an infinite series of payoffs results in the mutual defection, but (2) a pair of players learning from a certain finite series of payoffs results in mutual cooperation. 
Although we should be sufficiently careful to conclude the general outcomes in the iterated prisoner's dilemma with learners, 
these two cases indicates the possibility of a new type of mechanism leading to mutual cooperation which has not been 
reported in past studies \cite{Nowak2006}.
Namely, mutual cooperation in the iterated prisoner's dilemma can be achieved by learning with a modestly length-limited series of payoffs.

\section{Conclusion \label{sec-Conclusion}}
Analysis of a game of players with probabilistic strategies which depend on a long history requires analysis of a higher order Markov process. Often, however, we wish to understand the game using just a low-order marginal stationary distribution rather than the full-ordered stationary.
For such cases, namely determining a low-order marginal probability of a high-order Markov process, we explored here the nature of a series of Markov processes constructured in a recursive manner.
The main result, Theorem \ref{thm-RecursiveMarginal}, states the condition that we can analyze a low-order marginal stationary distribution of the corresponding infinite order Markov process in a closed form. 
The numerical studies in the previous section confirm the computational advantage of $m$-shift stationary distributions over the stationary distributions of the corresponding $k^{\text{th}}$-order Markov chains.
Our technique is potentially applicable to the analysis of a class of game with players with probabilistic strategies which depend on a series of past states.
We expect it will open a new field of study into types of games with complex players who learn from their past experience.

\section*{Acknowledgments}
This study was supported by JSPS KAKENHI 23300099, 15KT0013, 16H01609, and 16H05860.

\bibliographystyle{spmpsci}      
\bibliography{SimpleClass}

\section{Appendix: Approximation in the exchangeable special case \label{sec-Exchangeable}}

In the special case $\alpha_{i} = 1$ for each $i$ of Definition \ref{def-LearningGames}, we can closely approximate a stationary distribution for a relatively large $k$. This approximation was employed in the previous analysis \cite{Hidaka2015}.
Exploiting the exchangeability of actions in a state series, the size of the state space of this special case is $(k+1)k \ldots (k-N+3)$ for $N \ge 2$, which is quite smaller than the size of the original state space, $N^{k}$. With the following formulation, we can compute the stationary distribution for a relatively large $k$.

With $\alpha_{i} = 1$ for each $i$, the order of joint states in a series is inconsequential, because reward weights are equal at every step. In this case, we can identify two joint states 
$\left( X_{t}, X_{t-1}, \ldots, X_{t-s}, \ldots, X_{t-s'}, \ldots \right)$ and $\left( X_{t}, X_{t-1}, \ldots, X_{t-s'}, \ldots, X_{t-s}, \ldots \right)$ for any pair $s, s' < \infty$. Thus, for a finite $k$, we rewrite state space by the counts of $N$ joint states, \[
 C_{t-k}^{t} := \{ (C_{1}, C_{2}, \ldots, C_{N}) : C_{i} = 
| \{ x \ \in \{ X_{t-k}, X_{t-k+1}, \ldots, X_{t} \} : x = i  \} | \}.
\]
Over this counting state space $C_{t-k}^{t}$, 
we obtain the recursive equation on the stationary distribution
\begin{equation}
 P( C_{1}, C_{2}, \ldots, C_{N} ) = \sum_{i=1}^{N}
\pi_{i}P( C_{1} - \delta_{i,1}, C_{2} - \delta_{i,2}, \ldots, C_{N} - \delta_{i,N} ).
\label{eq-CountingRecursiveEquation}
\end{equation}
For $i = 1, \hdots, N$, $\pi_{i}$ is the conditional probability 
\[
 \pi_{i} = P( C_{1}, C_{2}, \ldots, C_{N} \mid  C_{1} - \delta_{i,1}, C_{2} - \delta_{i,2}, \ldots, C_{N} - \delta_{i,N} ).
\]
Given a probability of some initial state, we can compute the forward-in-time probabilities over these counting states using (\ref{eq-CountingRecursiveEquation})
until $k$ is sufficiently large. In the numerical implementation of this case, we removed counting states with a probability less than $10^{-10}$ for computational efficiency. This rounding reduced the probability by less than $1\%$ of the total probability 1 in our analysis.

\end{document}